\RequirePackage{ifpdf}
\ifpdf 
\documentclass[pdftex]{sigma}
\else
\documentclass{sigma}
\fi

\numberwithin{equation}{section}
\numberwithin{theorem}{section}
\numberwithin{proposition}{section}
\numberwithin{lemma}{section}
\numberwithin{corollary}{section}
\numberwithin{definition}{section}
\numberwithin{example}{section}
\numberwithin{remark}{section}
\numberwithin{note}{section}

\usepackage{mathrsfs}

\begin{document}

\allowdisplaybreaks

\renewcommand{\PaperNumber}{003}

\FirstPageHeading

\ShortArticleName{Intertwinors on Functions over the Product of Spheres}

\ArticleName{Intertwinors on Functions over the Product of Spheres}

\Author{Doojin HONG}

\AuthorNameForHeading{D. Hong}

\Address{Department of Mathematics, University of North Dakota, Grand Forks ND 58202, USA}  
\Email{\href{mailto:doojin.hong@und.edu}{doojin.hong@und.edu}}

\ArticleDates{Received August 23, 2010, in f\/inal form December 30, 2010;  Published online January 06, 2011}

\Abstract{We give explicit formulas for the intertwinors on the scalar functions over the
product of spheres with the natural pseudo-Riemannian product metric using the spectrum
generating technique. As a consequence, this provides another proof of the even
order conformally invariant dif\/ferential operator formulas obtained
earlier by T.~Branson and the present author.}


\Classification{53A30; 53C50}

\section{Introduction}

Branson, \'Olafsson, and {\O}rsted presented in \cite{BOO:96} a method of computing intertwining operators between principal series representations induced from maximal parabolic subgroups of semisimple Lie groups in the case where $K$-types occur with multiplicity at most one. One of the main ideas is that the intertwining relation, when ``compressed'' from a $K$-type to another $K$-type, can provide a purely numerical relation between eigenvalues on the $K$-types being considered through some
relatively simple calculations. This procedure of getting a recursive numerical
relation is referred to as ``spectrum generating'' technique.

When the group under consideration acts as conformal transformations, the intertwinors are conformally invariant operators (see \cite[Section~3]{BOO:96} for example) and they
have been one of the major subjects in mathematics and physics.

Explicit formulas for such operators on many manifolds are potentially important.
For instance, the precise form of
Polyakov formulas in even dimensions for the quotient of functional determinants
of operators only depends on some
constants that appear in the spectral asymptotics of the operators in question~\cite{Branson:95}.

In 1987, Branson~\cite{Branson:87} presented explicit formulas for invariant operators on functions and dif\/ferential
forms over the double cover $S^1\times S^{n-1}$ of the $n$-dimensional compactif\/ied
Minkowski space. And Branson and Hong \cite{BH:06, Hong:00, Hong:04} gave explicit determinant quotient formulas
for operators on spinors and twistors including the Dirac and Rarita--Schwinger operators over
$S^1\times S^{n-1}$.

In this paper, we show that the spectrum generating technique can be applied to get explicit formulas for the intertwinors on the scalar functions
over general product of spheres, $S^p\times S^q$ with the natural
pseudo-Riemannian metric.

\section{Some background on conformally covariant operators}

We brief\/ly review conformal covariance and the intertwining
relation (for more details, see~\cite{Branson:96, BOO:96}).

Let $(M,g)$ be an $n$-dimensional pseudo-Riemannian manifold.
If $f$ is a (possibly local) dif\/feomorphism on $M$, we denote by $f\cdot$ the
natural action of $f$ on tensor f\/ields which acts on vector f\/ields as
$f\cdot X=(df)X$ and on covariant tensors as $f\cdot \phi=(f^{-1})^*\phi$.

A vector f\/ield $T$ is said to be {\it conformal} with {\it conformal factor}
$\omega\in C^{\infty}(M)$
if
\begin{gather*}
\mathcal{L}_T g=2\omega g ,
\end{gather*}
where $\mathcal{L}$ is the Lie derivative. The conformal vector f\/ields form a
Lie algebra $\mathfrak{c}(M,g)$.
A~conformal transformation on $(M,g)$ is a (possibly local) dif\/feomorphism $h$
for which $h\cdot g=\Omega^2 g$ for some positive function $\Omega\in
C^{\infty}(M)$. The global conformal transformations form a group~$\mathscr{C}(M,g)$.
We have representations def\/ined by
\begin{gather*}
\mathfrak{c}(M,g)\stackrel{U_a}{\longrightarrow}  \mbox{End}\, C^{\infty}(M),\qquad
 U_a(T)=\mathcal{L}_T+a\omega\qquad \mbox{and}\\
\mathscr{C}(M,g)\stackrel{u_a}{\longrightarrow}  \mbox{Aut}\, C^{\infty}(M), \qquad
 u_a(h)=\Omega^a h\cdot
\end{gather*}
for $a\in \mathbb{C}$.

Note that if a conformal vector f\/ield $T$ integrates to a one-parameter group of
global conformal transformation $\{h_{\varepsilon}\}$, then
\begin{gather*}
\{U_a(T)\phi\}(x)=\frac{d}{d\varepsilon}\Big|_{\varepsilon=0}\{u_a(h_{-\varepsilon})\phi\}(x)  .
\end{gather*}
In this sense, $U_a$ is the inf\/initesimal representation corresponding to $u_a$.

A dif\/ferential operator $D:C^{\infty}(M)\rightarrow C^{\infty}(M)$ is said to be
{\it infinitesimally conformally covariant of bidegree} $(a,b)$ if
\begin{gather*}
DU_a(T)\phi=U_b(T)D\phi
\end{gather*}
for all $T\in \mathfrak{c}(M,g)$ and $D$ is said to be {\it conformally covariant of bidegree} $(a,b)$ if
\begin{gather*}
Du_a(h)\phi=u_b(h)D\phi
\end{gather*}
for all $h\in \mathscr{C}(M,g)$.

To relate conformal covariance to conformal invariance, we let $M$ be an $n$-dimensional manifold with metric $g$ and recall \cite{Branson:96, Branson:97} that
conformal weight of a subbundle $V$ of some tensor bundle over $M$ is
$r\in \mathbb{C}$
if and only if
\begin{gather*}
\bar{g}=\Omega^2 g \ \Longrightarrow \ \overline{g_{_V}}=\Omega^{-2r} g_{_V}  ,
\end{gather*}
where $\Omega > 0 \in C^\infty(M)$ and $g_{_V}$ is the induced bundle metric from the metric
$g$.
Tangent bundle, for instance, has conformal weight~$-1$. Let us denote a bundle $V$ with
conformal weight $r$ by~$V^r$. Then we can impose new conformal weight $s$ on~$V^r$ by
taking tensor product of it with the bundle $I^{(s-r)/n}$ of scalar $((s-r)/n)$-densities. Now if
we view an operator of bidegree~$(a,b)$ as an operator from the bundle with conformal weight~$-a$ to the bundle with conformal weight~$-b$, the operator becomes conformally
invariant.

As an example, let us consider the conformal Laplacian on $M$:
\begin{gather*}
Y=\triangle+\frac{n-2}{4(n-1)}R,
\end{gather*}
where $\triangle=-g^{ab}\nabla_a\nabla_b$ and $R$ is the scalar curvature.
Note that $Y: C^\infty(M) \rightarrow C^\infty(M)$ is conformally covariant of bidegree $((n-2)/2,(n+2)/2)$. That is,
\begin{gather*}
\overline{Y}=\Omega^{-\frac{n+2}{2}}Y\mu\big(\Omega^{\frac{n-2}{2}}\big)  ,
\end{gather*}
where $\overline{Y}$ is $Y$ evaluated in $\overline{g}$ and
$\mu\big(\Omega^{\frac{n-2}{2}}\big)$ is multiplication by $\Omega^{\frac{n-2}{2}}$.
If we let $V=C^\infty(M)$ and view $Y$ as an operator
\begin{gather*}
Y: \ V^{-\frac{n-2}{2}} \rightarrow V^{-\frac{n+2}{2}}  ,
\end{gather*}
we have, for $\phi\in V^{-\frac{n-2}{2}}$,
\begin{gather*}
\overline{Y}\,\overline{\phi}=\overline{Y\phi}  ,
\end{gather*}
where $\overline{Y}$, $\overline{\phi}$, and $\overline{Y\phi}$ are $Y$, $\phi$, and $Y\phi$ computed in $\overline{g}$, respectively.

\section{Conformal structure and intertwining relation}

Let us f\/irst brief\/ly review conformal structure on $S^p\times S^q$ (for more details, see \cite{TB:02} for example).

Let $\mathbb{R}^{p+1,q+1}$ be the $(n+2)$-dimensional pseudo-Riemannian manifold
$\mathbb{R}^{p+q+2}$ equipped with the pseudo-Riemannian metric
\begin{gather*}
-\xi_{-p}^2-\cdots-\xi_0^2+\xi_1^2+\cdots+\xi_{q+1}^2
\end{gather*}
and def\/ine submanifolds of $\mathbb{R}^{p+1,q+1}$ by
\begin{gather*}
\Xi:=\{-\xi_{-p}^2-\cdots-\xi_0^2+\xi_1^2+\cdots+\xi_{q+1}^2=0\}\setminus \{0\}  ,\\
M:=\{\xi_{-p}^2+\cdots+\xi_0^2=\xi_1^2+\cdots+\xi_{q+1}^2=1\}\simeq S^p\times S^q  .
\end{gather*}
The natural action of the multiplicative group $\mathbb{R}^\times_+=\{r\in\mathbb{R}: r>0\}$ on~$\Xi$ identif\/ies $\Xi/\mathbb{R}^\times_+$ with~$M$ and
the natural action of the orthogonal group $G:=O(p+1,q+1)$ on  $\mathbb{R}^{p+1,q+1}$ stabilizes the metric cone~$\Xi$ and these two actions commute. So we get a $G$-equivariant principal $\Xi/\mathbb{R}^\times_+$ bundle:
\begin{gather*}
\Phi: \ \Xi\rightarrow M, \qquad (\xi_{-p},\dots,\xi_0,\xi_1,\dots,\xi_{q+1})\rightarrow
\frac{1}{\sqrt[]{\xi_{-p}^2+\cdots+\xi_0^2}}(\xi_{-p},\dots,\xi_0,\xi_1,\dots,\xi_{q+1})  .
\end{gather*}
The standard pseudo-Riemannian metric
on $\mathbb{R}^{p+1,q+1}$ induces the standard pseudo-Riemannian metric
$-g_{{}_{S^p}}+g_{{}_{S^q}}$ on $S^p\times S^q\simeq M$.
For $h\in G$, $z\in M$, the map
\begin{gather*}
z\stackrel{h}{\longrightarrow} h\cdot z\stackrel{\Phi}{\longrightarrow} \Phi(h\cdot z)
\end{gather*}
is a conformal transformation on $M$ and $G$ is the full conformal group of $S^p\times S^q$.

A basis of the conformal vector f\/ields on $S^p\times S^q$ in homogeneous coordinates is
\begin{gather*}
L_{\alpha\beta}=\varepsilon_\alpha\xi_\alpha\partial_\beta-\varepsilon_\beta\xi_\beta\partial_\alpha   ,
\end{gather*}
where $\partial_\alpha=\partial / \partial\xi_\alpha$ and $-\varepsilon_p=\cdots =-\varepsilon_0=\varepsilon_1=\cdots =\varepsilon_{q+1}=1$.

The subalgebra spanned by $L_{\alpha\beta}$ for $-p\le \alpha,\beta \le 0$ is a copy of $\mathfrak{so}(p+1)$ so it generates $\operatorname{SO}(p+1)$ group of isometries. Likewise,
the $L_{\alpha\beta}$ for $1\le \alpha, \beta =q+1$ generate $\operatorname{SO}(q+1)$ group of isometries.

Let $\operatorname{SO}_0(p+1,q+1)$ be the identity component of $G$ and note $K=\operatorname{SO}(p+1)\times\operatorname{SO}(q+1)$ is
a~maximal compact subgroup of $\operatorname{SO}_0(p+1,q+1)$. Elements in $K$ act as isometries on
$S^p\times S^q$ and
proper conformal vector f\/ields are the ones with mixed indices:
$L_{\alpha\beta}$ for $-p\le \alpha\le 0 < \beta \le q$.

To express the typical proper conformal vector f\/ield $L_{01}$ in intrinsic coordinates
of $S^p\times S^q$, we
let $\tau$ be the azimuthal angle in $S^p$. That is, set
\begin{gather*}
\xi_0=\cos \tau, \qquad 0 \le \tau \le \pi,
\end{gather*}
and complete $\tau$ to a set of spherical angular coordinates $(\tau, \tau_1,\dots,\tau_{p-1})$ on $S^p$.
Likewise, set
\begin{gather*}
\xi_1=\cos \rho, \qquad 0 \le \rho \le \pi,
\end{gather*}
on $S^q$ and complete $\rho$ to a set of spherical angular coordinates $(\rho, \rho_1, \dots, \rho_{q-1})$ on $S^q$.
Then
\begin{gather}
L_{01}=\cos\rho\sin \tau\partial_\tau+\cos \tau\sin\rho\partial_\rho =:T,\nonumber\\
\omega_{01}=\cos \tau\cos\rho=:\varpi .\label{Tw}
\end{gather}
Note that $\sin\tau\partial_\tau$ (resp. $\sin\rho\partial_\rho$) is conformal on the Riemannian
$S^p$ (resp.~$S^q$) with the conformal factor $\cos\tau$ (resp.~$\cos\rho$).

Let $A=A_{2r}$ be an intertwinor of order $2r$ of the $(\mathfrak{g},K)$ representation on functions over $S^p\times S^q$.
That is, $A$ is a $K$-map satisfying
the intertwining relation~\cite{Branson:96,BOO:96}
\begin{gather}\label{rel-1}
A\left(\mathcal{L}_X+\left(\frac{n}2-r\right)\varpi\right)=
\left(\mathcal{L}_X+\left(\frac{n}2+r\right)\varpi\right)A\qquad \text{for all} \ \ X\in \mathfrak{g}  ,
\end{gather}
where $\mathcal{L}_X$ is the {\em Lie derivative}.

\begin{remark}
If the intertwinor acts on tensors of $\left(\begin{array}{c}l\\m\end{array}\right)$-type, then the Lie derivative should be changed to $\mathcal{L}_X+(l-m)\varpi$ in (\ref{rel-1}) \cite[p.~347]{Branson:97}. Note also that
$\mathcal{L}_X=\nabla_X$ on functions.
\end{remark}

The spectrum generating relation that converts (\ref{rel-1}) is given in the following lemma.
\begin{lemma}
Let $T$ and $\varpi$ be as in \eqref{Tw}.
Consider the Riemannian--Bochner Laplacian $\nabla^{*,{\rm R}}\nabla:=-g^{\alpha\beta}\nabla_\alpha\nabla_\beta=:N$, where $g=g_{{}_{S^p}}+g_{{}_{S^q}}$ is the standard Riemannian metric on $S^p\times S^q$. Then,
\begin{gather}
[N,\varpi]=2\left(\nabla_T+\frac{n}2\varpi\right)\nonumber
\end{gather}
on tensors of any type. Here $[\cdot,\cdot]$ is the usual operator commutator.
\end{lemma}

\begin{proof}
If $\varphi$ is any smooth section, then
\begin{gather*}
[N,\varpi]\varphi = (\triangle \varpi)\varphi-2\iota (d\varpi) \nabla\varphi
 = \cos\rho\triangle_{S^p}\cos\tau+\cos\tau \triangle_{S^q}\cos\rho+2\nabla_T  \varphi\\
\phantom{[N,\varpi]\varphi}{}
= (p+q)\varpi+2\nabla_T  \varphi  ,
\end{gather*}
where $\iota$ is the interior multiplication and both $\triangle_{S^p}$ and $\triangle_{S^q}$
are Riemannian Laplacians.
\end{proof}

Thus the intertwining relation (\ref{rel-1}) becomes
\begin{gather}\label{rel-2}
A\left(\frac{1}{2}[N,\varpi]-r\varpi\right)=
\left(\frac{1}{2}[N,\varpi]+r\varpi\right)A  .
\end{gather}

 Recall that the space of $j$-th order spherical harmonics on the Riemannian $S^p$ is the irreducible $\operatorname{SO}(p+1)$-module def\/ined by
\begin{gather*}
E(j)=\{f\in C^{\infty}(S^p):\triangle_{S^p}f=j(j+p-1)f\}
\end{gather*}
and the space $L^2(S^p)$ decomposes as
\begin{gather*}
L^2(S^p)\simeq\bigoplus_{j=0}^{\infty}E(j) .
\end{gather*}
Let $F(k)$ be the space of $k$-th order spherical harmonics on the Riemannian $S^q$ and def\/ine
\begin{gather*}
\mathcal{V}(j,k):=E(j)\otimes F(k) .
\end{gather*}
Note that we have a multiplicity free $K$-type decomposition into
$K$-f\/inite subspaces of the space of smooth functions on $S^p\times S^q$:
\begin{gather*}
 \bigoplus_{j,k=0}^{\infty}\mathcal{V}(j,k)
\end{gather*}
and $K$ operator $A$ acts as a scalar multiplication on each $\mathcal{V}(j,k)$.

On $S^p$, a proper conformal factor maps $E(j)$ to the
direct sum $E(j+1)\oplus E(j-1)$. See~\cite{Branson:97} for details.
Thus, a proper conformal factor on $S^p\times S^q$ maps a $K$-type $\mathcal{V}(j,k)$ to land in the direct sum of
4 types $\mathcal{V}(j',k')$:
\begin{gather}\label{arrows}
\left(\begin{array}{ll}
\mathcal{V}(j-1,k+1)\quad & \mathcal{V}(j+1,k+1) \\
\mathcal{V}(j-1,k-1) & \mathcal{V}(j+1,k-1)
\end{array}\right)  .
\end{gather}
Let $\alpha$ be the $K$-type $\mathcal{V}(j,k)$ and $\beta$ be any $K$-type appearing in (\ref{arrows}). We apply the intertwining relation (\ref{rel-2}) to a section $\varphi$ in $\alpha$:
\begin{gather*}
A\left(\frac{1}{2}[N,\varpi]-r\varpi\right)\varphi=
\left(\frac{1}{2}[N,\varpi]+r\varpi\right)A\varphi \\
\qquad{} \Leftrightarrow \
A\left(\frac{1}{2}(N(\varpi\varphi)-\varpi(N_\alpha\varphi))-r\varpi\varphi\right)=\mu_\alpha
\left(\frac{1}{2}(N(\varpi \varphi)-\varpi(N_\alpha\varphi))+r\varpi \varphi\right)  ,
\end{gather*}
where $\mu_\alpha$ (resp. $N_\alpha$) is the eigenvalue of $A$ (resp. $N$) on the $K$-type $\alpha$.
Note that $\varpi\varphi$ is a~direct sum of the $K$-types appearing in
(\ref{arrows}).
Let Proj${}_\beta\varpi|_a\varphi$ be the projection of $\varpi\varphi$ onto the $K$-type~$\beta$.
The ``compression'', from the $K$-type $\alpha$ to the $K$-type $\beta$,
of the above relation becomes \underline{Proj${}_\beta\varpi|_\alpha$ times}
\begin{gather}\label{numer}
\left(\frac12N|^\beta_\alpha+r\right)\mu_\alpha
=\left(\frac12N|^\beta_\alpha-r\right)\mu_\beta ,
\end{gather}
where $\mu_\beta$ (resp.~$N_\beta$) is the eigenvalue of $A$ (resp.~$N$) on
the $K$-type $\beta$ and $N|^\beta_\alpha:=N_\beta-N_\alpha$.  The underlined phrase above is a key point. We
have achieved a factorization in which one factor is purely numerical
(that appearing in~(\ref{numer})).  ``Canceling'' the other factor,
Proj${}_\beta\varpi|_\alpha$, we get purely numerical recursions that
are guaranteed to give intertwinors.  If we wish to see the {\em uniqueness}
of intertwinors this way, we need to establish the nontriviality of
the Proj${}_\beta\varpi|_\alpha$. In fact this nontriviality follows from
Branson \cite[Section~6]{Branson:97}.

To compute $N|^\beta_\alpha$, we need the following lemma.
\begin{lemma}
Let $\alpha=\mathcal{V}(j,k)$ and $\beta=\mathcal{V}(j',k')$.
\begin{gather*}
N|^\beta_\alpha=j'(p-1+j')+k'(q-1+k')-j(p-1+j)-k(q-1+k)  .
\end{gather*}
\end{lemma}

\begin{proof}
On the Riemannian $S^p$, $\nabla^*\nabla$ acts as (see~\cite{Branson:92} for details)
\begin{gather*}
j(p-1+j)\quad  \text{on} \ \ E(j)\qquad \text{and} \qquad
N|_{\mathcal{V}(j,k)}=\nabla^*\nabla|_{E(j)}+\nabla^*\nabla|_{F(k)}  .\tag*{\qed}
\end{gather*}
\renewcommand{\qed}{}
\end{proof}

Let $J=j+\frac{p-1}{2}$ and $K=k+\frac{q-1}{2}$ for $j,k\in \mathbb{N}$. Then
the {\it transition quantities} $\mu_\beta/ \mu_a$ are
\begin{gather*}
\left(\begin{array}{ll}
\dfrac{-J+K+1+r}{-J+K+1-r}\quad & \dfrac{J+K+1+r}{J+K+1-r}\vspace{2mm}\\
\dfrac{-J-K+1+r}{-J-K+1-r}& \dfrac{J-K+1+r}{J-K+1-r}
\end{array}\right)
\end{gather*}
relative to (\ref{arrows}).

Note that $\mathcal{V}(f',j')$ can be reached from $\mathcal{V}(f,j)$ if and only
if $|f'-f|=|j'-j|=1$. So~$\mathcal{E}^0$, the direct sum of $\mathcal{V}(j,k)$ with
$j+k$ even, is a $(\mathfrak{g}, K)$-invariant subspace of $\mathcal{E}(K,\lambda_0)$,
as is the corresponding odd space $\mathcal{E}^1$.
Choosing normalization $\mu_{{}_{\mathcal{V}(0,0)}}=1$ (resp.~$\mu_{{}_{\mathcal{V}(1,0)}}=1$) on $\mathcal{E}^0$ (resp.~$\mathcal{E}^1$), we get

\begin{theorem}
The unique spectral function $Z_\varepsilon(r;f,j)$ on $\mathcal{E}^\varepsilon$ is, up to
normalization,
\begin{gather*}
\dfrac{\Gamma(\frac{1}{2}(K+J+1+r))\Gamma(\frac{1}{2}(K-J+1+r))
\Gamma(\frac{1}{2}(\varepsilon-\frac{p-q}{2}+1-r))\Gamma(\frac{1}{2}(\varepsilon+\frac{p+q}{2}-r))}
{\Gamma(\frac{1}{2}(K+J+1-r))\Gamma(\frac{1}{2}(K-J+1-r))
\Gamma(\frac{1}{2}(\varepsilon-\frac{p-q}{2}+1+r))\Gamma(\frac{1}{2}(\varepsilon+\frac{p+q}{2}+r))}
  .
\end{gather*}
\end{theorem}

When $r$ is a positive integer, the spectral function provides conformally covariant
dif\/ferential operators. To see this, let
\begin{gather*}
B:=\sqrt{\triangle_{S^p}+\left(\dfrac{p-1}2\right)^2},\qquad
C:=\sqrt{\triangle_{S^q}+\left(\dfrac{q-1}2\right)^2},
\end{gather*}
so that $B$ and $C$ are nonnegative operators with
\begin{gather*}
\triangle_{S^p}=B^2-\left(\frac{p-1}2\right)^2,\qquad
\triangle_{S^q}=C^2-\left(\frac{q-1}2\right)^2.
\end{gather*}
The eigenvalue list for $\triangle_{S^p}$ \cite{Branson:92,IT:78} is
\begin{gather*}
j(p-1+j),\qquad j=0,1,2,\dots,
\end{gather*}
so the eigenvalue list for $B$ is
\begin{gather*}
j+\dfrac{p-1}2,\qquad j=0,1,2,\dots.
\end{gather*}
Similarly, the eigenvalue list for $C$ is
\begin{gather*}
k+\dfrac{q-1}2,\qquad k=0,1,2,\dots.
\end{gather*}
Note that for $r=1$, $Z_\varepsilon(1;f,j)$ is up to a constant
\begin{gather*}
\frac{1}{2}(K+J)\frac{1}{2}(K-J).
\end{gather*}
Thus $Z_\varepsilon(1;f,j)$ agrees with a constant multiple of the conformal Laplacian on
$S^p\times S^q$
\begin{gather*}
(C-B)(C+B) =\triangle_{S^q}
-\triangle_{S^p} +\left(\dfrac{q-1}2\right)^2-\left(\dfrac{p-1}2\right)^2
 =\triangle_{S^q}-\triangle_{S^p}+\dfrac{n-2}{4(n-1)}\text{Scal}  ,
\end{gather*}
where Scal = scalar curvature of $S^p\times S^q$.

In general, we have

\begin{corollary}[\cite{BH:07}]
For a positive integer $r$, $Z_\varepsilon(r;f,j)$ is a constant multiple of
the differential operator
\begin{gather*}
(C+B-r+1)\cdots (C+B+r-1)\cdot (C-B-r+1)\cdots (C-B+r-1) .
\end{gather*}
where the increments are by $2$ units each time.
\end{corollary}

\begin{remark}
The same technique can be applied to the case of dif\/ferential form bundles
over $S^p\times S^q$. But here certain $K$-types occur with multiplicity~2.
On the multiplicity~1 types, exactly the same method yields the spectral function.
\end{remark}

\pdfbookmark[1]{References}{ref}
\LastPageEnding


\begin{thebibliography}{99}

\footnotesize\itemsep=0pt

\bibitem{Branson:87}
Branson T.,
Group representations arising from Lorentz conformal geometry,
\href{http://dx.doi.org/10.1016/0022-1236(87)90025-5}{{\it J.~Funct. Anal.}}  {\bf 74} (1987), 199--291.

\bibitem{Branson:92}
Branson T.,
Harmonic analysis in vector bundles associated to the rotation and spin groups,
\href{http://dx.doi.org/10.1016/0022-1236(92)90050-S}{{\it J.~Funct. Anal.}} {\bf 106} (1992), 314--328.

\bibitem{Branson:95}
Branson T.,
Sharp inequalities, the functional determinant, and the complementary series,
\href{http://dx.doi.org/10.2307/2155203}{{\it Trans. Amer. Math. Soc.}} {\bf 347} (1995), 3671--3742.

\bibitem{Branson:96}
Branson T.,
Nonlinear phenomena in the spectral theory of geometric linear dif\/ferential operators,
in Quantization, Nonlinear Partial Dif\/ferential Equations, and Operator Algebra (Cambridge, MA, 1994), {\it  Proc. Sympos. Pure Math.}, Vol.~59, Amer. Math. Soc., Providence, RI, 1996, 27--65.

\bibitem{Branson:97}
Branson T.,
Stein--Weiss operators and ellipticity,
\href{http://dx.doi.org/10.1006/jfan.1997.3162}{{\it J. Funct. Anal.}} {\bf 151} (1997), 334--383.

\bibitem{BH:06}
Branson T.P., Hong D.,
Spectrum generating on twistor bundle,
{\it Arch. Math. (Brno)} {\bf 42} (2006), suppl., 169--183,
\href{http://arxiv.org/abs/math.DG/0606524}{math.DG/0606524}.

\bibitem{BH:07}
Branson T., Hong D.,
Translation to bundle operators,
\href{http://dx.doi.org/10.3842/SIGMA.2007.102}{{\it SIGMA}} {\bf 3} (2007), 102, 14~pages,
\href{http://arxiv.org/abs/math.DG/0606552}{math.DG/0606552}.

\bibitem{BOO:96}
Branson T., \'Olafsson G.,  {\O}rsted B.,
Spectrum generating operators, and intertwining operators for representations induced from a maximal parabolic subgroups,
\href{http://dx.doi.org/10.1006/jfan.1996.0008}{{\it J.~Funct. Anal.}} {\bf 135} (1996), 163--205.

\bibitem{Hong:00}
Hong D.,
Eigenvalues of Dirac and Rarita--Schwinger operators,
in Clif\/ford Algebras (Cookeville, TN, 2002), {\it Prog. Math. Phys.}, Vol.~34, Birkh\"auser Boston, Boston, MA, 2004, 201--210.

\bibitem{Hong:04}
Hong D.,
Spectra of higher spin operators,
Ph.D. Thesis, University of Iowa, 2004.

\bibitem{IT:78}
Ikeda A., Taniguchi Y.,
Spectra and eigenforms of the Laplacian on $S^n$ and $P^n(\mathbb{C})$,
{\it Osaka~J. Math.} {\bf 15} (1978), 515--546.

\bibitem{TB:02}
Kobayashi T., {\O}rsted B.,
Analysis on the minimal representation of $O(p,q)$. I.~Realization via conformal geometry,
\href{http://dx.doi.org/10.1016/S0001-8708(03)00012-4}{{\it Adv. Math.}} {\bf 180} (2003), 486--512,
\href{http://arxiv.org/abs/math.RT/0111083}{math.RT/0111083}.

\end{thebibliography}
\end{document}